\documentclass[11pt]{amsart}

\usepackage{amsmath,amssymb,amsthm}
\usepackage{graphicx}
\usepackage{hyperref}
\usepackage{geometry}
\usepackage{amssymb}
\usepackage{amsmath}
\usepackage{amsfonts}
\usepackage{hyperref}
\usepackage{graphicx}
\usepackage{float}
\usepackage{amsmath,cite}
\usepackage{bookmark}
\usepackage{xcolor}
\usepackage{placeins}
\usepackage{subcaption} 
\usepackage{multirow}
\usepackage{array}
\usepackage{enumitem}
\usepackage{textcomp}
\usepackage{wasysym}
\usepackage{geometry}
\usepackage{comment}

\title{Performance comparison of Python, MATLAB and R for numerical solutions of SI and SIR epidemiological models}
\title[Computational Comparison of Numerical Methods for SI/SIR Models]{Performance comparison of Python, MATLAB and R for numerical solutions of SI and SIR epidemiological models}
\author{Berkay Özışık}
\address{Ankara University, Faculty of Sciences, Department of Mathematics, Besevler, 06100, Ankara, TÜRKİYE}
\email{berkayozisik0@gmail.com}

\author{Elif Demirci}
\address{Ankara University, Faculty of Sciences, Department of Mathematics, Besevler, 06100, Ankara, TÜRKİYE}
\email{edemirci@ankara.edu.tr}

\begin{document}
\maketitle

\begin{abstract}
Mathematical modeling plays a vital role in epidemiology, offering insights into the spread and control of infectious diseases. The compartmental models developed by Kermack and McKendrick, particularly the SI (Susceptible–Infected) and SIR (Susceptible–Infected–Recovered) models, form the basis of many epidemic studies. While some simple cases permit analytical solutions, most real-world models require numerical methods such as Euler’s method, the fourth-order Runge-Kutta (RK4) method, and Predictor–Corrector (P–C) methods. These methods are typically implemented in scientific computing software like Python, MATLAB, and R. However, the computational efficiency and run-time performance of these software tools in solving epidemiological models have not been comprehensively compared in the literature. This study addresses this gap by solving the SI and SIR models using Euler’s method, RK4, and P–C methods in Python, MATLAB, and R. Execution times are recorded for each implementation to evaluate computational efficiency. Additionally, for the SI model—where an exact analytical solution exists—$R^2$ values are computed to assess numerical accuracy. For the SIR model, a high-accuracy reference solution is obtained by solving the system using MATLAB’s ODE45 solver, and the SIR solutions computed via the RK4 method in MATLAB are compared against this reference. The results provide a comparative perspective on the accuracy and run-time performance across different software and numerical methods, offering practical guidance for researchers and practitioners in selecting suitable tools for epidemic modeling.
\end{abstract}

\section{Introduction}
Mathematical models, which represent real-world problems, are frequently used in science and engineering. Epidemiology is an important field where mathematical modeling plays a crucial role, particularly in understanding and controlling the spread of infectious diseases. A foundational approach in epidemic modeling is the compartmental model developed by Kermack and McKendrick in 1927 (Kermack and McKendrick, 1927), which remains influential today. Two basic epidemiological models, namely SI and SIR models, divide the population into two (Susceptible–Infected) and three (Susceptible–Infected–Recovered) compartments, respectively. Although some simple epidemic models allow us to obtain their exact solutions using elementary methods, the more complex models require numerical methods. Euler’s method, the fourth-order Runge-Kutta method (RK4) and predictor-corrector (P-C) methods are common numerical methods used to solve SI and SIR models numerically. Applying these numerical methods manually is evidently impractical due to their computational complexity. To apply numerical methods, various software has been used to facilitate such calculations. Among these software, Python, MATLAB, and R are widely used in mathematical and scientific computing. However, when utilizing numerical methods in these software, run-time performance becomes a critical factor. Several studies have investigated the numerical solutions of SIR epidemic models using various methods and software. In (Side et al., 2018), MATLAB was used to apply the RK4 method to a SIR model with tuberculosis transmission without providing any information on computational time. Moreover, in (Iskandar et al., 2022) a SIR model was considered for COVID-19 parameters in Malaysia and the model was solved with fourth order Runge-Kutta method in MATLAB, numerically and run-time of the program was also not reported. However, (Ashgi et al., 2021) presented a SIR model using COVID-19 data from Wuhan, China. This model with given parameters was solved with both Euler’s method and the RK4 method in Python, numerically. This study includes the comparison of execution times of the methods. The literature review showed that there is no study where the SI and SIR epidemiological models were solved using the Euler’s method, the RK4 method, and the P-C method in all three software: Python, MATLAB, and R, simultaneously. In addition, no study was found that compares the run-time of the software by applying these numerical methods in the existing literature. This paper aims to compare the run-time performance of three prominent programming languages—Python, MATLAB, and R—when employed for the numerical solutions of SI and SIR models. We use Euler’s method, RK4 method and P-C method in each software to solve the models. We also analyse the $R^2$ values for the SI model, that we could find the exact solution, for each numerical method, and for the SIR model, the solution obtained via MATLAB’s ODE45 solver is taken as a reference, and the SIR results computed using the RK4 method in MATLAB are assessed relative to this reference.

\section{Methodology}
In this paper, we discuss two classical epidemiological models (SI and SIR models), formulated with initial conditions and parameter values as outlined below.\\

We consider the Susceptible-Infected (SI) model given by the following system of differential equations
\begin{equation}
\left\{
\begin{aligned}
\frac{dS}{dt} &= -\alpha SI =: f(t,S,I), \\
\frac{dI}{dt} &= \alpha SI =: g(t,S,I)
\end{aligned}
\right.
\qquad t \in I := [t_0,b]
\end{equation}
with initial conditions
\[
S(t_0) = S_0, \quad I(t_0) = I_0, \quad N(t_0) = N_0.
\]

Here, \( S \) and \( I \) represent the number of susceptible and infected individuals in the population, respectively. The total population is \( N = S + I \). The parameter \( \alpha \) denotes the transmission parameter of the infection.
\\
The analytical solution of system (1) is: 
\[
\begin{aligned}
S(t) &= \frac{N \, e^{-\alpha t N} \, c}{1 + e^{-\alpha t N} \, c}, \\
I(t) &= \frac{N(N - S_0)}{N - S_0 + e^{-N \alpha t} S_0} 
\quad \text{where } 
c = \frac{S_0}{N - S_0}
\end{aligned}
\]

\vspace{0.1em}
The Susceptible-Infected-Recovered (SIR) model given by the following system of differential equations
\begin{equation}
\left\{
\begin{aligned}
\frac{dS}{dt} &= -\alpha SI =: f(t,S,I), \\
\frac{dI}{dt} &= \alpha SI - \beta I =: g(t,S,I), \\
\frac{dR}{dt} &= \beta I =: u(t,I)
\end{aligned}
\right.
\qquad t \in I := [t_0,b]
\end{equation}
with initial conditions:
\[
S(t_0) = S_0, \quad I(t_0) = I_0, \quad R(t_0) = R_0, \quad N(t_0) = N_0.
\]

Here, \( S \), \( I \), and \( R \) denote the susceptible, infected, and recovered populations, respectively. The total population is \( N = S + I + R \). The parameter \( \alpha \) is the transmission parameter, and \( \beta \) is the recovery parameter.

\subsection{Numerical Methods}
In order to solve the aforementioned models given by (1) and (2) we apply Euler’s method, the fourth-order Runge-Kutta (RK4) method, and Predictor-Corrector (P-C) method to the initial value problems (Chapra and Canale, 2015).

\vspace{0.5em}

Consider the following initial value problem:
\begin{equation}
\frac{dy}{dx} = f(x, y), \quad y(x_0) = y_0, \quad x \in I := [x_0, b].
\end{equation}
\\[1ex]
\noindent
Euler's method is a simple, first order numerical  method that the solution of an initial value problem given as (3) is evaluated iteratively as: 
\begin{equation}
y_{n+1} = y_n + h f(x_n, y_n), \quad n \geq 0
\end{equation}

Here, \( n \) represents the number of steps and \( h \) represents the step size, defined by:
\begin{equation}
h := \frac{b - x_0}{n}.
\end{equation}

\vspace{0.1em}

The formulation obtained by applying  the Euler's method to the models given in (1) and (2) are 
\vspace{0.1em}
\vspace{0.1em}
\vspace{0.1em}

\begin{equation}
\left\{
\begin{aligned}
S_{n+1} &= S_n + h f(t_n,S_n,I_n)=S_n+h(-\alpha S_nI_n), \\
I_{n+1} &= I_n + h g(t_n,S_n,I_n)=I_n+h(\alpha S_nI_n) \\
\end{aligned}
\right.
\end{equation}

and

\begin{equation}
\left\{
\begin{aligned}
S_{n+1} &= S_n + h f(t_n,S_n,I_n)=S_n+h(-\alpha S_nI_n), \\
I_{n+1} &= I_n + h g(t_n,S_n,I_n)=I_n+h(\alpha S_nI_n-\beta I_n), \\
R_{n+1} &= R_n + h u(t_n,I_n)=R_n+h(\beta I_n),
\end{aligned}
\right.
\end{equation} 
respectively.

\noindent
\vspace{0.5em}

The RK4 method is a fourth-order numerical method for finding numerical solutions of IVP's and the RK4 formulation is (Chapra and Canale, 2015):

\[
y_{n+1} = y_n + \frac{h}{6} (k_1+2 k_2+2k_3+k_4), \quad n \geq 0
\]
where, 
\[
k_1:=f(x_n,y_n),
\]
\[
k_2:=f(x_n+\frac{h}{2},y_n+h\frac{k_1}{2}),
\]
\[
k_3:=f(x_n+\frac{h}{2},y_n+h\frac{k_2}{2}),
\]
\[
k_4:=f(x_n+h,y_n+hk_3).
\]

\noindent

Thus, RK4 method can be formulated for the models (1) and (2) by:

\begin{equation}
\left\{
\begin{aligned}
S_{n+1} &= S_n + \frac{h}{6} (k_1^{S}+2k_2^{S}+ 2k_3^{S}+k_4^{S}), \\
I_{n+1} &= I_n + \frac{h}{6} (k_1^{I}+2k_2^{I}+ 2k_3^{I}+k_4^{I}), \\
\end{aligned}
\right.
\end{equation}
where,
{\scriptsize
\begin{align*}
k_1^{S} &:= f(t_n, S_n, I_n)=-\alpha S_nI_n, \\
k_2^{S} &:= f\left(t_n + \frac{h}{2}, S_n+\frac{hk_1^{S}}{2}, I_n+\frac{hk_1^{I}}{2}\right)=-\alpha(S_n + \frac{hk_1^S}{2})(I_n + \frac{hk_1^I}{2}), \\
k_3^{S} &:= f\left(t_n + \frac{h}{2}, S_n+\frac{hk_2^{S}}{2}, I_n+\frac{hk_2^{I}}{2}\right)=-\alpha(S_n + \frac{hk_2^S}{2})(I_n + \frac{hk_2^I}{2}), \\
k_4^{S} &:= f(t_n + h, S_n+hk_3^{S} ,I_n+hk_3^{I})=-\alpha(S_n + hk_3^S)(I_n + hk_3^I), \\
k_1^{I} &:= g(t_n, S_n, I_n)=\alpha S_nI_n, \\
k_2^{I} &:= g\left(t_n + \frac{h}{2}, S_n+\frac{hk_1^{S}}{2}, I_n+\frac{hk_1^{I}}{2}\right)=\alpha(S_n + \frac{hk_1^S}{2})(I_n + \frac{hk_1^I}{2}), \\
k_3^{I} &:= g\left(t_n + \frac{h}{2}, S_n+\frac{hk_2^{S}}{2}, I_n+\frac{hk_2^{I}}{2}\right)=\alpha(S_n + \frac{hk_2^S}{2})(I_n + \frac{hk_2^I}{2}), \\
k_4^{I} &:= g(t_n + h, S_n+hk_3^{S} ,I_n+hk_3^{I})=\alpha(S_n + hk_3^S)(I_n + hk_3^I)
\end{align*}
}

and

\begin{equation}
\left\{
\begin{aligned}
S_{n+1} &= S_n + \frac{h}{6} (k_1^{S}+2k_2^{S}+ 2k_3^{S}+k_4^{S}), \\
I_{n+1} &= I_n + \frac{h}{6} (k_1^{I}+2k_2^{I}+ 2k_3^{I}+k_4^{I}), \\
R_{n+1} &= R_n + \frac{h}{6} (k_1^{R}+2k_2^{R}+ 2k_3^{R}+k_4^{R}), \\
\end{aligned}
\right.
\end{equation}
where,
{\scriptsize
\begin{align*}
& k_1^S := f(t_n, S_n, I_n) = -\alpha S_nI_n, \\
& k_2^S := f(t_n + \frac{h}{2}, S_n + \frac{hk_1^S}{2}, I_n + \frac{hk_1^I}{2}) = -\alpha(S_n + \frac{hk_1^S}{2})(I_n + \frac{hk_1^I}{2}), \\
& k_3^S := f(t_n + \frac{h}{2}, S_n + \frac{hk_2^S}{2}, I_n + \frac{hk_2^I}{2}) = -\alpha(S_n + \frac{hk_2^S}{2})(I_n + \frac{hk_2^I}{2}), \\
& k_4^S := f(t_n + h, S_n+hk_3^{S} ,I_n+hk_3^{I})=-\alpha(S_n + hk_3^S)(I_n + hk_3^I),\\
& k_1^I := g(t_n, S_n, I_n) = \alpha S_nI_n -\beta I_n, \\
& k_2^I := g(t_n + \frac{h}{2}, S_n + \frac{hk_1^S}{2}, I_n + \frac{hk_1^I}{2}) = \alpha(S_n + \frac{hk_1^S}{2})(I_n + \frac{hk_1^I}{2})- \beta (I_n+\frac{hk_1^I}{2}), \\
& k_3^I := g(t_n + \frac{h}{2}, S_n + \frac{hk_2^S}{2}, I_n + \frac{hk_2^I}{2})=\alpha(S_n + \frac{hk_2^S}{2})(I_n + \frac{hk_2^I}{2})- \beta (I_n+\frac{hk_2^I}{2}), \\
& k_4^I := g(t_n + h, S_n+hk_3^{S} ,I_n+hk_3^{I})=\alpha(S_n + hk_3^S)(I_n + hk_3^I)- \beta (I_n+hk_3^I), \\
& k_1^R := u(t_n, I_n) = \beta I_n, \\
& k_2^R := u(t_n + \frac{h}{2}, I_n + \frac{hk_1^I}{2}) = \beta (I_n+ \frac{hk_1^I}{2}), \\
& k_3^R := u(t_n + \frac{h}{2}, I_n + \frac{hk_2^I}{2}) = \beta (I_n+ \frac{hk_2^I}{2}), \\
& k_4^R := u(t_n + h, I_n + hk_3^I) = \beta (I_n+ hk_3^I).\\
\end{align*}
}

The P-C method is an iterative technique where a predictor is first used to estimate the next value, and then a corrector refines this estimate. Predictor and corrector are defined as:
\begin{align*}
\textbf{Predictor:} \quad & \tilde{y}_{i+1} = y_i + h f(x_i, y_i), \\
\textbf{Corrector:} \quad & y_{i+1} = y_i + \frac{h}{2} \left[ f(x_i, y_i) + f(x_{i+1}, \tilde{y}_{i+1}) \right].
\end{align*}
where $\tilde{y}_{i+1}$ denotes the predicted value obtained from the predictor step, which is subsequently refined by the corrector (Chapra and Canale, 2015).

If we extend the P-C Method to (1) and (2), we get:
\begin{equation}
\left\{
\begin{aligned}
S_{n+1} &= S_n + \frac{h}{2} \left[ f(t_n, S_n, I_n) + f(t_{n+1}, \tilde{S}_{n+1}, \tilde{I}_{n+1}) \right], \\
I_{n+1} &= I_n + \frac{h}{2} \left[ g(t_n, S_n, I_n) + g(t_{n+1}, \tilde{S}_{n+1}, \tilde{I}_{n+1}) \right]
\end{aligned}
\right.
\end{equation}

and
\vspace{0.5em}

\begin{equation}
\left\{
\begin{aligned}
S_{n+1} &= S_n + \frac{h}{2} \left[ f(t_n, S_n, I_n) + f(t_{n+1}, \tilde{S}_{n+1}, \tilde{I}_{n+1}) \right], \\
I_{n+1} &= I_n + \frac{h}{2} \left[ g(t_n, S_n, I_n) + g(t_{n+1}, \tilde{S}_{n+1}, \tilde{I}_{n+1}) \right], \\
R_{n+1} &= R_n + \frac{h}{2} \left[ u(t_n, I_n) + u(t_{n+1}, \tilde{I}_{n+1}) \right],
\end{aligned}
\right.
\end{equation}
where $\tilde{S}_{n+1}=S_n+hf(t_n,S_n,I_n)$, $\tilde{I}_{n+1}=I_n+hg(t_n,S_n,I_n)$.

\subsection{$R^2$ Formula}
In order to measure the errors between the numerical solutions and the exact solution we use the $R^2$ formula which is outlined in the formula (12).
\begin{equation}
R^2
= 1 - \frac{\sum_{i=1}^{n} \left( y_i - \tilde{y}_i \right)^2}
{\sum_{i=1}^{n} \left( y_i - \bar{y} \right)^2},
\end{equation}
where $y_i$ denotes the exact (reference) solution values, $\tilde{y}_i$ denotes the predicted values obtained by the numerical method, and $\bar{y}$ is the mean of the exact (reference) solution.

\subsection{Implementation and Time Measurement}
The execution times of the algorithms were measured by placing a timing function at the beginning and at the end of each implementation. Only the pure computational time of the numerical algorithms was recorded; operations such as defining parameters, initializing variables, model setup, and plotting were not included in the measurements. All simulations were performed on a MacBook Air equipped with an Apple M4 chip and 16 GB of RAM.

\section{Main Results}
In this study, we applied three numerical methods, Euler’s Method, RK4 (Runge-Kutta 4th order method), and the Predictor-Corrector (P-C) method, to solve the SI model (as given by equation 1) using three different software platforms: Python, MATLAB, and R.  We use the parameter value $\alpha = 2.18 \times 10^{-3}$ and initial conditions $S(0)=762$, $I(0)=1$, for $I:=[0,14]$ (days) from (Murray, 2002) which belongs to (Communicable Disease Surveillance Centre et al., 1978).

Since (1) can be solved analytically, we provide a comparison of the $R^2$ values for three different methods in three different software. The graphs of the numerical solutions of the SI model in software Python, MATLAB and R are shown in Figures 1, 2 and 3, respectively.
The $R^2$ values for different step sizes for each numerical methods in three different software are given in Table 1. Run-time of the programs for each numerical methods in three software are recorded as in Table 2.
Initially, we present the graphs of the numerical solutions of the SI model. 

\begin{figure}[H]
    \centering
    \begin{subfigure}{0.3\textwidth}
        \centering
        \includegraphics[width=\linewidth]{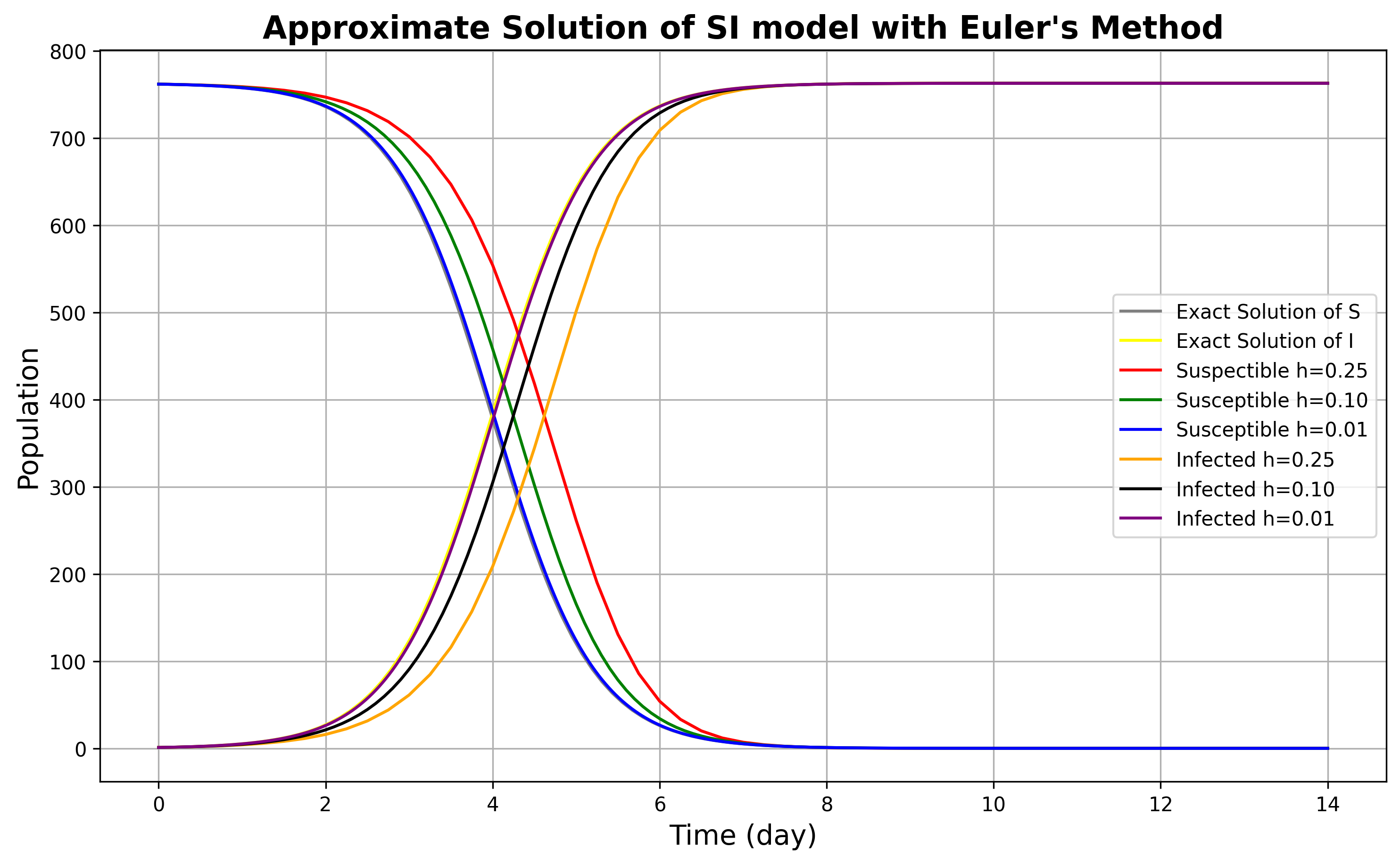}
        \caption{Solution of SI model with Euler's Method}
        \label{fig:euler_graph}
    \end{subfigure}
    \hfill 
    \begin{subfigure}{0.3\textwidth}
        \centering
        \includegraphics[width=\linewidth]{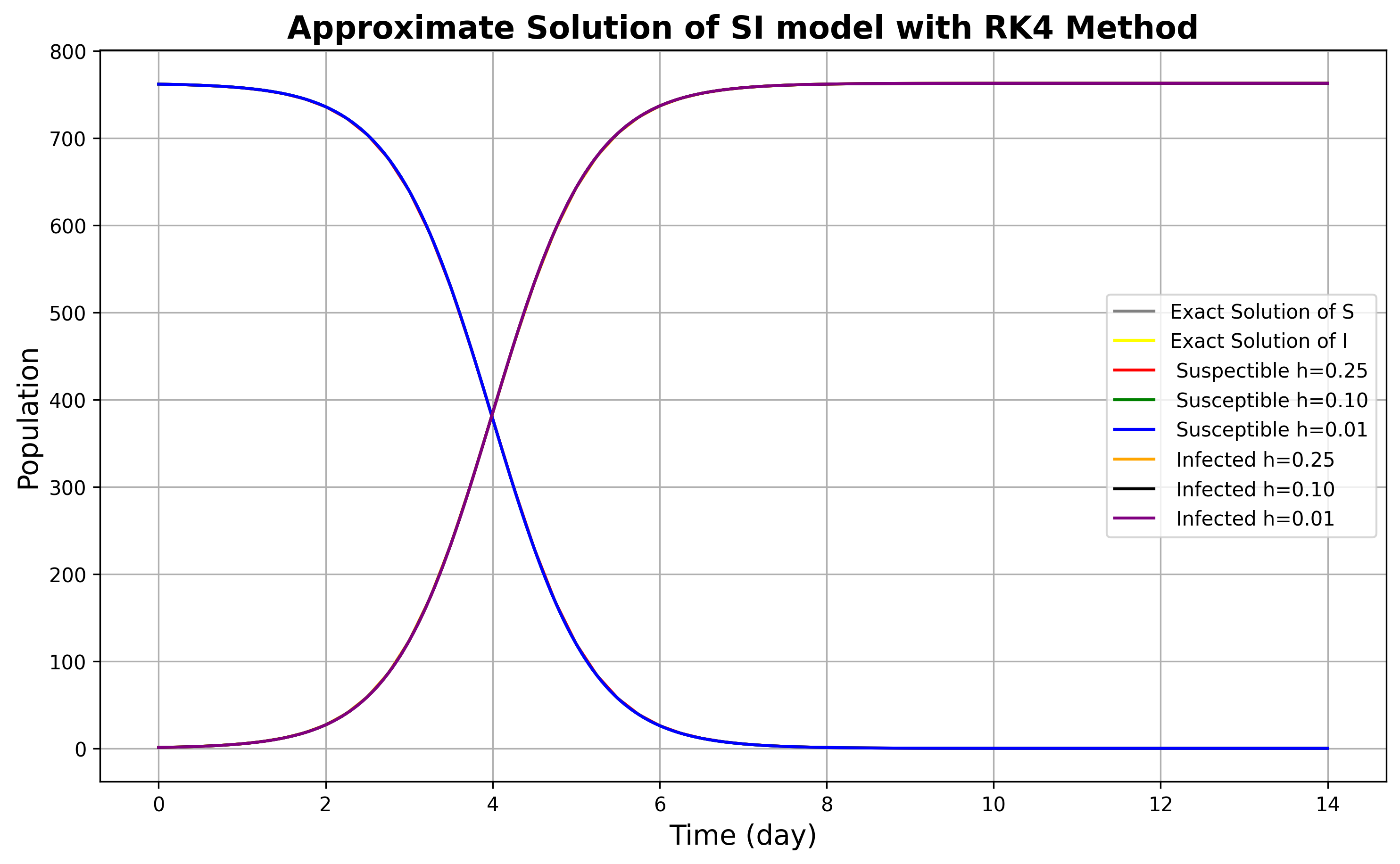}
        \caption{Solution of SI model with RK4 Method}
        \label{fig:RK4_graph}
    \end{subfigure}
    \hfill 
    \begin{subfigure}{0.3\textwidth}
        \centering
        \includegraphics[width=\linewidth]{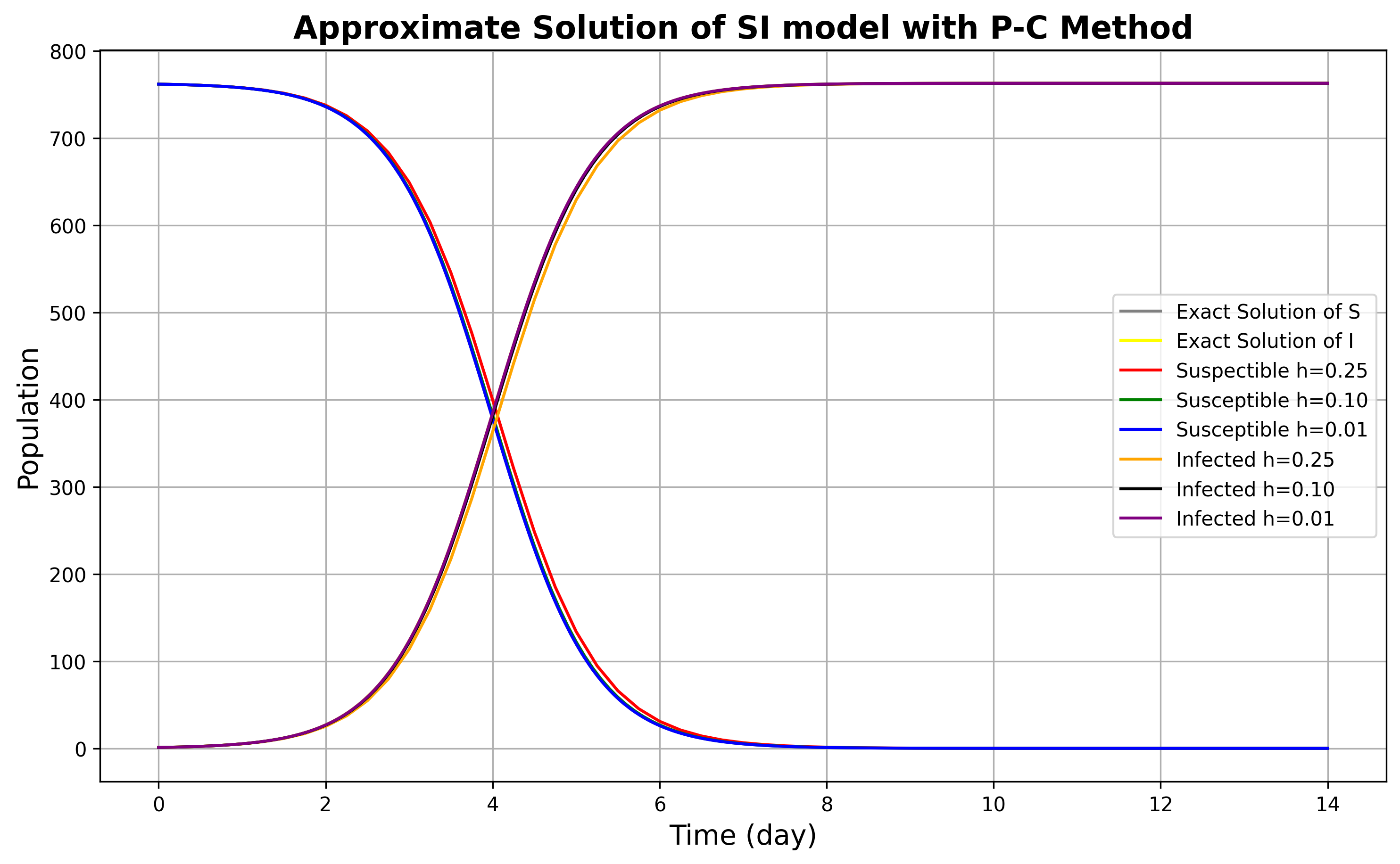}
        \caption{Solution of SI model with P-C Method}
        \label{fig:Adams_Graph}
    \end{subfigure}
    \caption{The graphs of the exact solution and the numerical solutions of the SI model in Python, using Euler's method, RK4 method and the P-C method.}
    \label{fig:all_graphs}
\end{figure}

\begin{figure}[H]
    \centering
    \begin{subfigure}{0.3\textwidth}
        \centering
        \includegraphics[width=\linewidth]{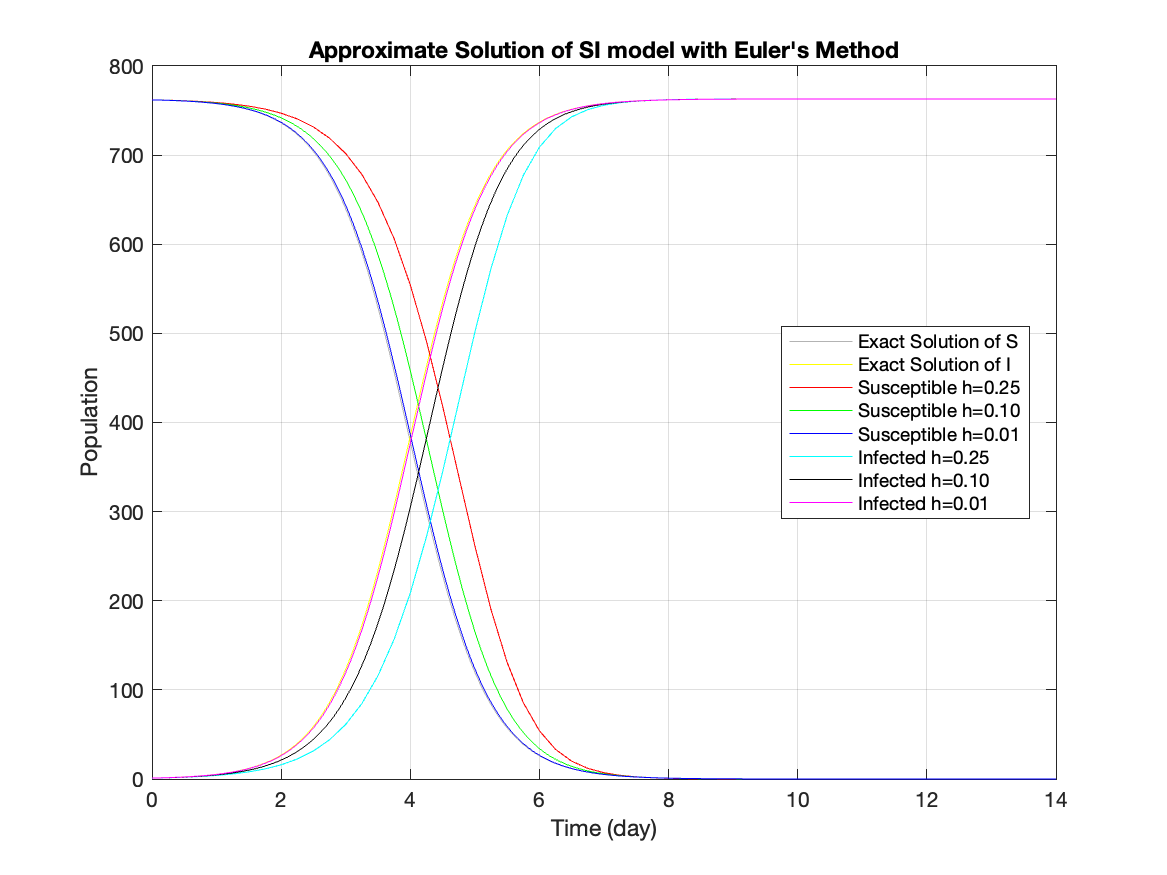}
        \caption{Solution of SI model with Euler's Method}
        \label{fig:m_euler_graph}
    \end{subfigure}
    \hfill 
    \begin{subfigure}{0.3\textwidth}
        \centering
        \includegraphics[width=\linewidth]{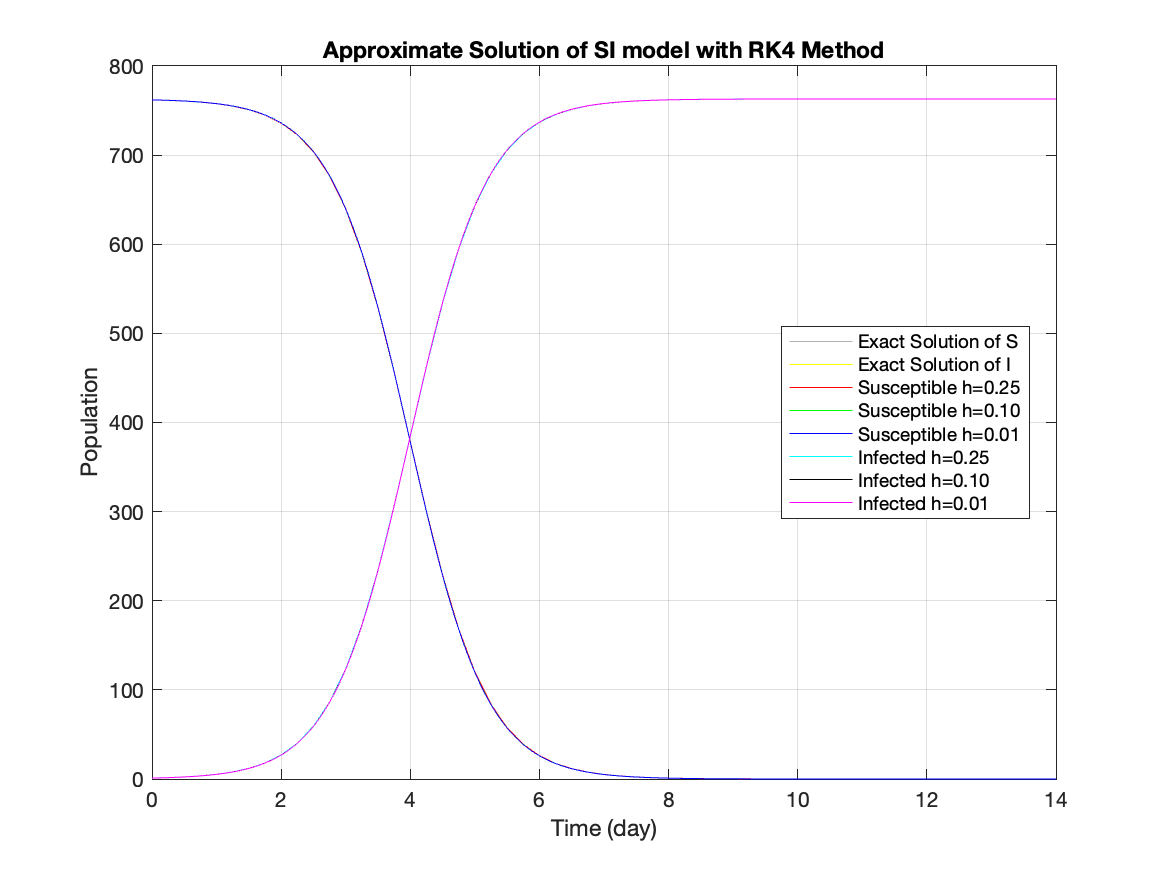}
        \caption{Solution of SI model with RK4 Method}
        \label{fig:m_RK4_graph}
    \end{subfigure}
    \hfill 
    \begin{subfigure}{0.3\textwidth}
        \centering
        \includegraphics[width=\linewidth]{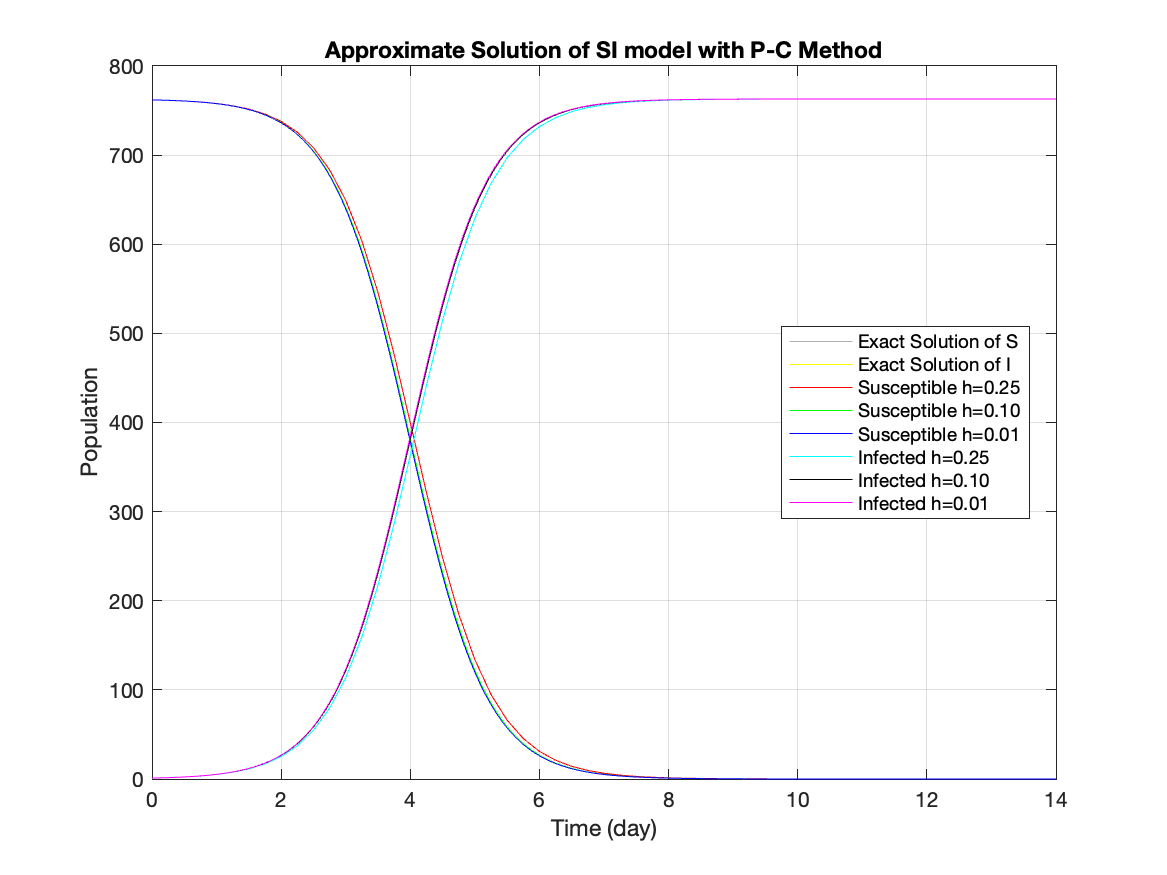}
        \caption{Solution of SI model with P-C Method}
        \label{fig:m_Adams_Graph}
    \end{subfigure}
    \caption{The graphs of the exact solution and the numerical solutions of the SI model in MATLAB, using Euler's method, RK4 method and the P-C method.}
    \label{fig:all_graphs_2}
\end{figure}

\begin{figure}[H]
    \centering
    \begin{subfigure}{0.3\textwidth}
        \centering
        \includegraphics[width=\linewidth]{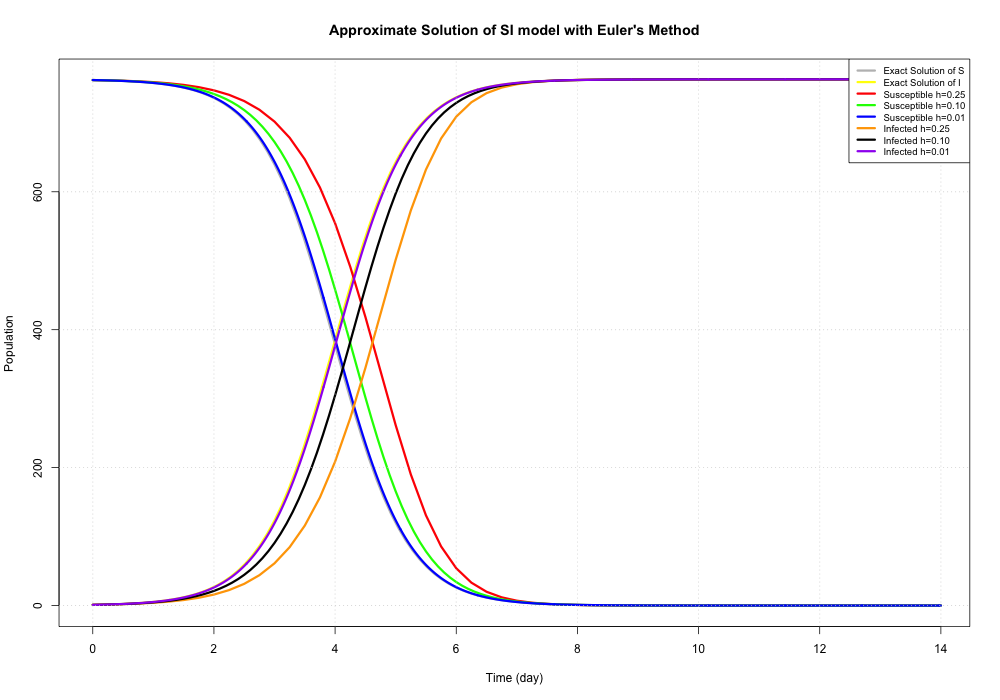}
        \caption{Solution of SI model with Euler Method}
        \label{fig:r_euler_graph}
    \end{subfigure}
    \hfill 
    \begin{subfigure}{0.3\textwidth}
        \centering
        \includegraphics[width=\linewidth]{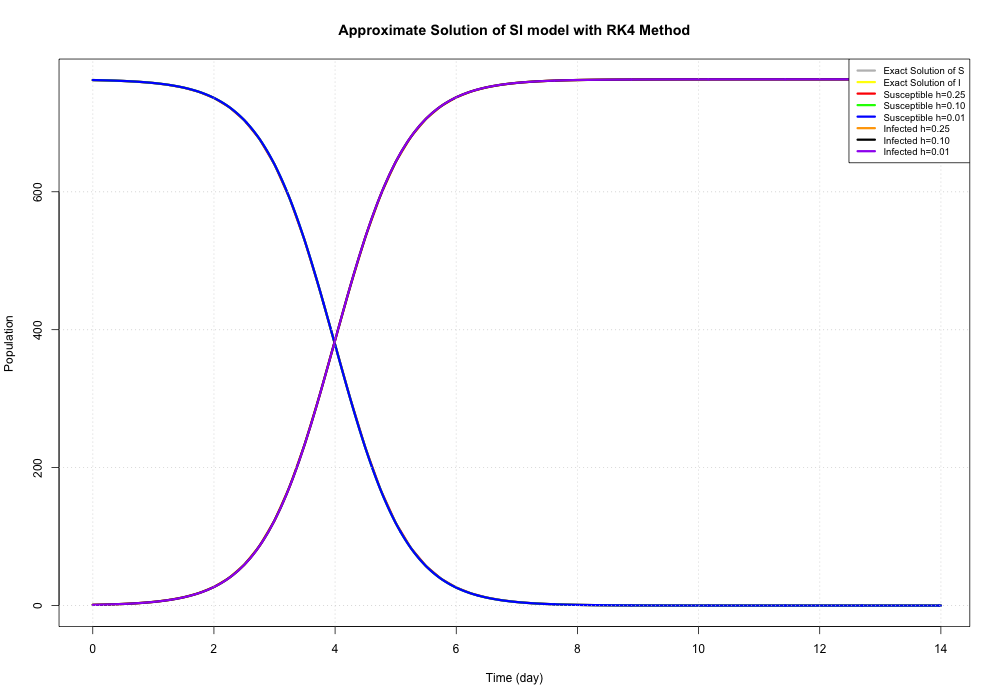}
        \caption{Solution of SI model with RK4 Method}
        \label{fig:r_RK4_graph}
    \end{subfigure}
    \hfill 
    \begin{subfigure}{0.3\textwidth}
        \centering
        \includegraphics[width=\linewidth]{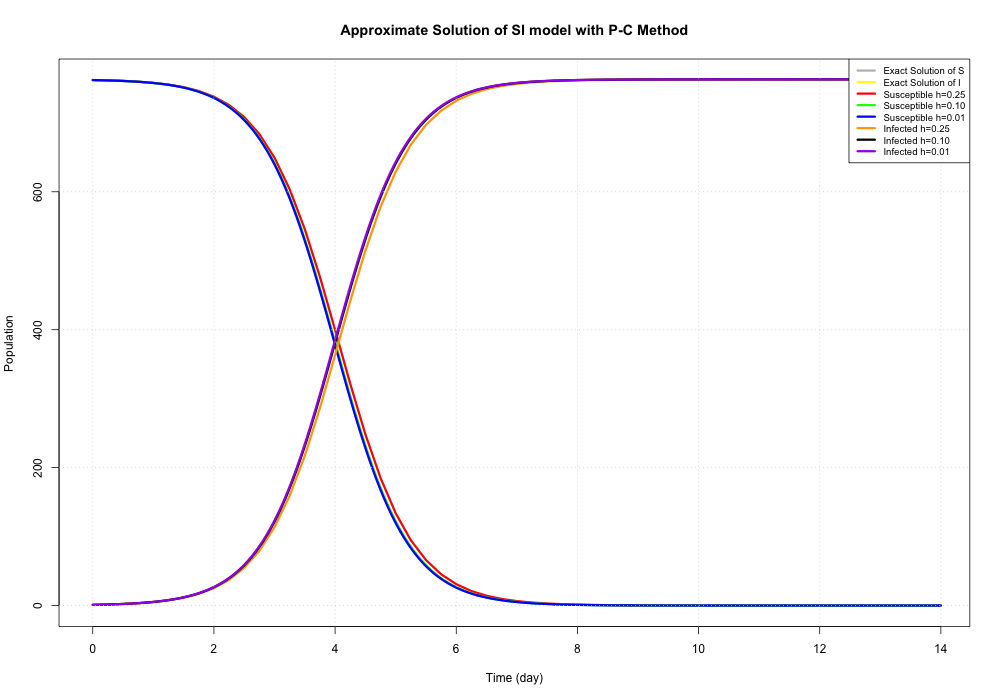}
        \caption{Solution of SI model with P-C Method}
        \label{fig:r_Adams_Graph}
    \end{subfigure}
    \caption{The graphs of the exact solution and the numerical solutions of the SI model in R, using Euler's method, RK4 method and the P-C method.}
    \label{fig:all_graphs}
\end{figure}

Now, we present the errors of the numerical solutions by using the $R^2$ measure formula.

\begin{table}[H]
\centering
\begin{tabular}{|>{\centering\arraybackslash}m{1.5cm}|c|c|c|c|}
    \hline
    \multirow{7}{*}{\rotatebox{90}{\textbf{Euler Method}}} &  & \textbf{Python} & \textbf{MATLAB} & \textbf{R} \\ \cline{2-5}
    & S &h=0.25: $R^2$ = 0.9585463 &h=0.25: $R^2$ = 0.9585463 &h=0.25: $R^2$ = 0.9585463 \\
    &   &h=0.10: $R^2$ = 0.9927564 &h=0.10: $R^2$ = 0.9927564 &h=0.10: $R^2$ = 0.9927564 \\
    &   &h=0.01: $R^2$ = 0.9999239 &h=0.01: $R^2$ = 0.9999239 &h=0.01: $R^2$ = 0.9999239 \\ \cline{2-5}
    & I &h=0.25: $R^2$ = 0.9585463 &h=0.25: $R^2$ = 0.9585463 &h=0.25: $R^2$ = 0.9585463 \\
    &   &h=0.10: $R^2$ = 0.9927564 &h=0.10: $R^2$ = 0.9927564 &h=0.10: $R^2$ = 0.9927564 \\
    &   &h=0.01: $R^2$ = 0.9999239 &h=0.01: $R^2$ = 0.9999239 &h=0.01: $R^2$ = 0.9999239 \\ \hline

    \multirow{7}{*}{\rotatebox{90}{\textbf{RK4 Method}}} &  & \textbf{Python} & \textbf{MATLAB} & \textbf{R} \\ \cline{2-5}
    & S &h=0.25: $R^2$ = 1.0 &h=0.25: $R^2$ = 1.0 &h=0.25: $R^2$ = 1.0 \\
    &   &h=0.10: $R^2$ = 1.0 &h=0.10: $R^2$ = 1.0 &h=0.10: $R^2$ = 1.0 \\
    &   &h=0.01: $R^2$ = 1.0 &h=0.01: $R^2$ = 1.0 &h=0.01: $R^2$ = 1.0 \\ \cline{2-5}
    & I &h=0.25: $R^2$ = 1.0 &h=0.25: $R^2$ = 1.0 &h=0.25: $R^2$ = 1.0 \\
    &   &h=0.10: $R^2$ = 1.0 &h=0.10: $R^2$ = 1.0 &h=0.10: $R^2$ = 1.0 \\
    &   &h=0.01: $R^2$ = 1.0 &h=0.01: $R^2$ = 1.0 &h=0.01: $R^2$ = 1.0 \\ \hline

    \multirow{7}{*}{\rotatebox{90}{\textbf{P - C Method}}} &  & \textbf{Python} & \textbf{MATLAB} & \textbf{R} \\ \cline{2-5}
    & S &h=0.25: $R^2$ = 0.9994189 &h=0.25: $R^2$ = 0.9994189 &h=0.25: $R^2$ = 0.9994189 \\
    &   &h=0.10: $R^2$ = 0.9999798 &h=0.10: $R^2$ = 0.9999798 &h=0.10: $R^2$ = 0.9999798 \\
    &   &h=0.01: $R^2$ = 1.0 &h=0.01: $R^2$ = 1.0 &h=0.01: $R^2$ = 1.0 \\ \cline{2-5}
    & I &h=0.25: $R^2$ = 0.9994189 &h=0.25: $R^2$ = 0.9994189 &h=0.25: $R^2$ = 0.9994189 \\
    &   &h=0.10: $R^2$ = 0.9999798 &h=0.10: $R^2$ = 0.9999798 &h=0.10: $R^2$ = 0.9999798 \\
    &   &h=0.01: $R^2$ = 1.0 &h=0.01: $R^2$ = 1.0 &h=0.01: $R^2$ = 1.0 \\ \hline
\end{tabular}
\caption{R$^2$ values of each numerical method in different software for the corresponding step sizes ($h$ values) for the SI model. (Values are rounded to seven decimal places.)}
\label{tab:numerical-errors}
\end{table}

Finally, we present the run-times of the programs that solves the SI model numerically.
\begin{table}[H]
\begin{center}
\begin{tabular}{|l|l|l|l|}
\hline
\textbf{} & \textbf{Python} & \textbf{MATLAB} & \textbf{R} \\
\hline
\textbf{Euler Method} 
& \begin{tabular}[c]{@{}l@{}}h=0.25: Time = 0.000037 \\h=0.10: Time = 0.000081\\h=0.01: Time = 0.000764
\end{tabular}
& \begin{tabular}[c]{@{}l@{}}h=0.25: Time = 0.004027
\\h=0.10: Time = 
0.007183
\\h=0.01: Time = 0.018744
\end{tabular}
& \begin{tabular}[c]{@{}l@{}}h=0.25: Time = 0.015112\\h=0.10: Time = 0.011526\\h=0.01: Time = 0.009426
\end{tabular} \\
\hline
\textbf{RK4 Method} 
& \begin{tabular}[c]{@{}l@{}}h=0.25: Time = 0.000163\\h=0.10: Time = 0.000398\\h=0.01: Time = 0.003950\end{tabular}
& \begin{tabular}[c]{@{}l@{}}h=0.25: Time = 0.012201
\\h=0.10: Time = 0.021967
\\h=0.01: Time = 0.025777
\end{tabular}
& \begin{tabular}[c]{@{}l@{}}h=0.25: Time = 0.024517\\h=0.10: Time = 0.024647\\h=0.01: Time = 0.031478\end{tabular} \\
\hline
\textbf{P - C Method} 
& \begin{tabular}[c]{@{}l@{}}h=0.25: Time = 0.000097\\h=0.10: Time = 0.000187\\h=0.01: Time = 0.001742\end{tabular}
& \begin{tabular}[c]{@{}l@{}}h=0.25: Time = 0.007261
\\h=0.10: Time = 0.015156
\\h=0.01: Time = 0.018358
\end{tabular}
& \begin{tabular}[c]{@{}l@{}}h=0.25: Time = 0.016191\\h=0.10: Time = 0.010234\\h=0.01: Time = 0.022058\end{tabular}
\\
\hline
\end{tabular}
\caption{Run-time of the software for each numerical method (seconds). (Values are rounded to six decimal places.)}
\end{center}
\end{table}

We now apply the numerical methods to the SIR model given by (2) by using the software Python, MATLAB and R. We use the parameter values $\alpha = 2.18 \times 10^{-3}$ and $\beta = (2.18 \times 10^{-3}) \times 202
$ as the transmission and recovery parameters, respectively and the initial values, $S(0)=762$, $I(0)=1$ and $R(0)=0$, for $I:=[0,14]$ (days) given in (Murray, 2002) that were obtained using the influenza epidemic data for a boys’ boarding school as reported in (Communicable Disease Surveillance Centre et al., 1978). \\[1ex]

The errors of the RK4 method for the SIR model with respect to the ODE45 reference solution are given in Table 3, while the graphs of the numerical solutions in Python, MATLAB, and R are presented in Figures 4, 5, and 6, respectively, and the run-time required for different step sizes is reported in Table 4.

\begin{figure}[H]
    \centering
    \begin{subfigure}{0.3\textwidth}
        \centering
        \includegraphics[width=\linewidth]{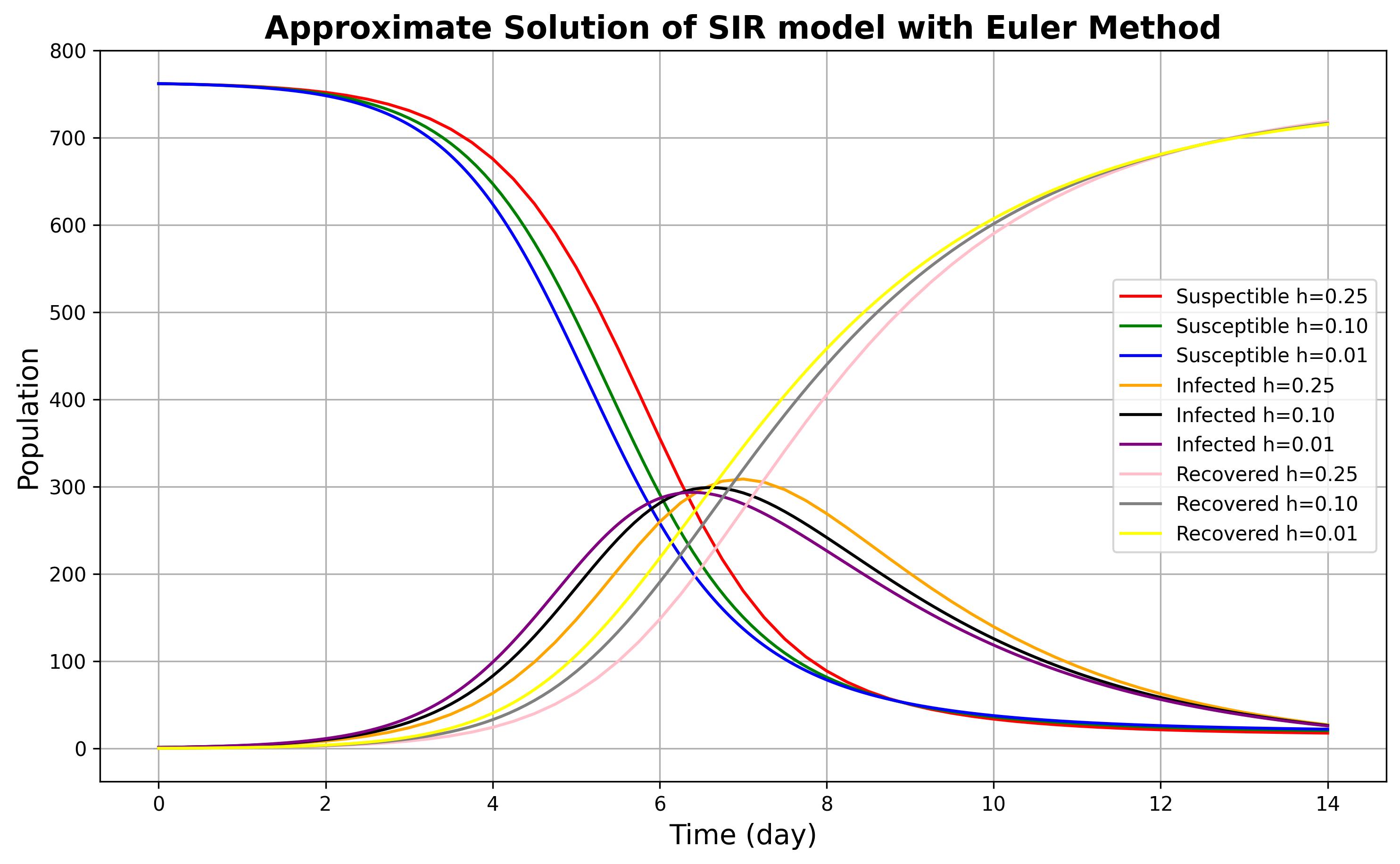}
        \caption{Solution of SIR model with Euler's Method}
        \label{fig:euler_graph2}
    \end{subfigure}
    \hfill 
    \begin{subfigure}{0.3\textwidth}
        \centering
        \includegraphics[width=\linewidth]{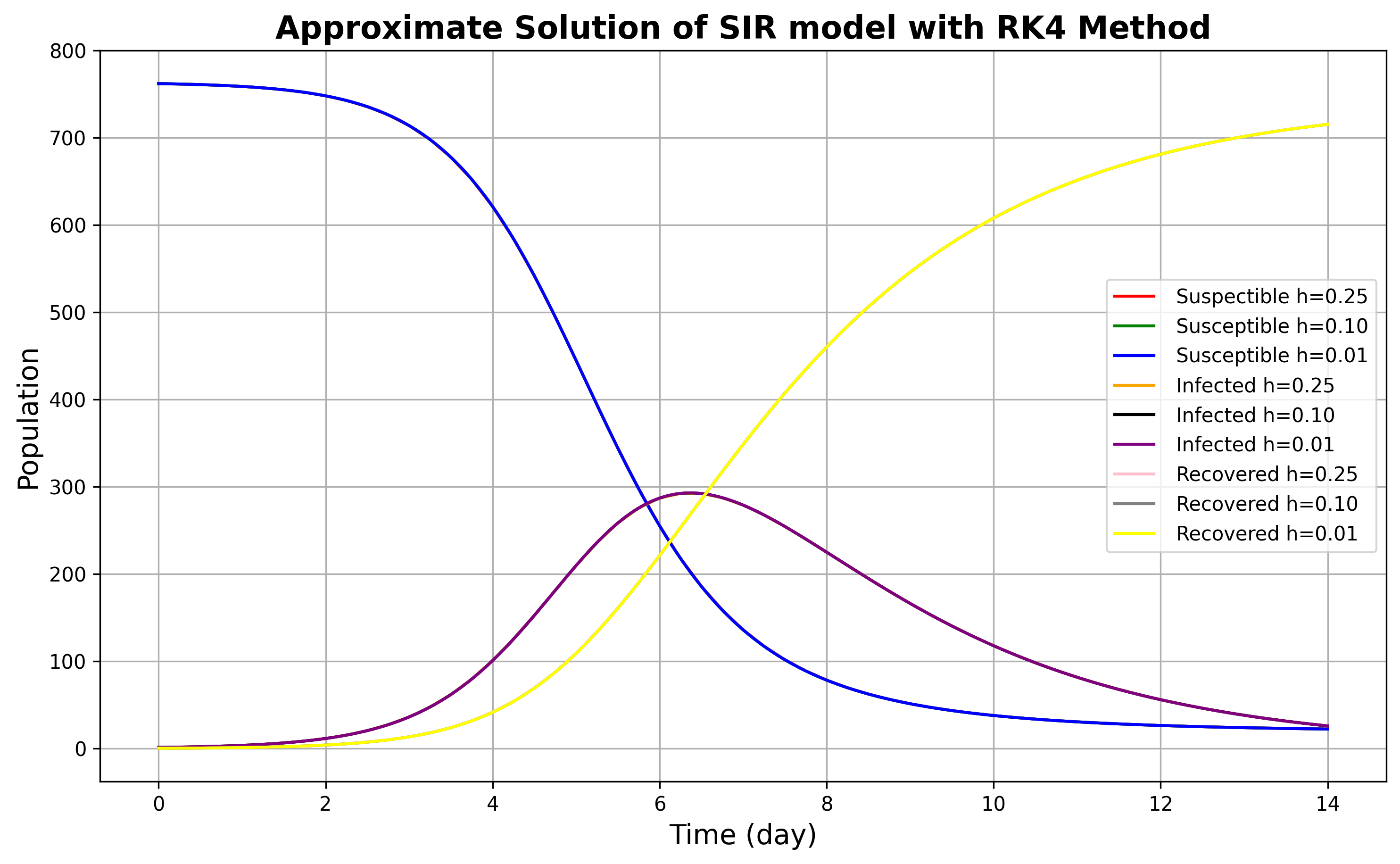}
        \caption{Solution of SIR model with RK4 Method}
        \label{fig:RK4_graph2}
    \end{subfigure}
    \hfill 
    \begin{subfigure}{0.3\textwidth}
        \centering
        \includegraphics[width=\linewidth]{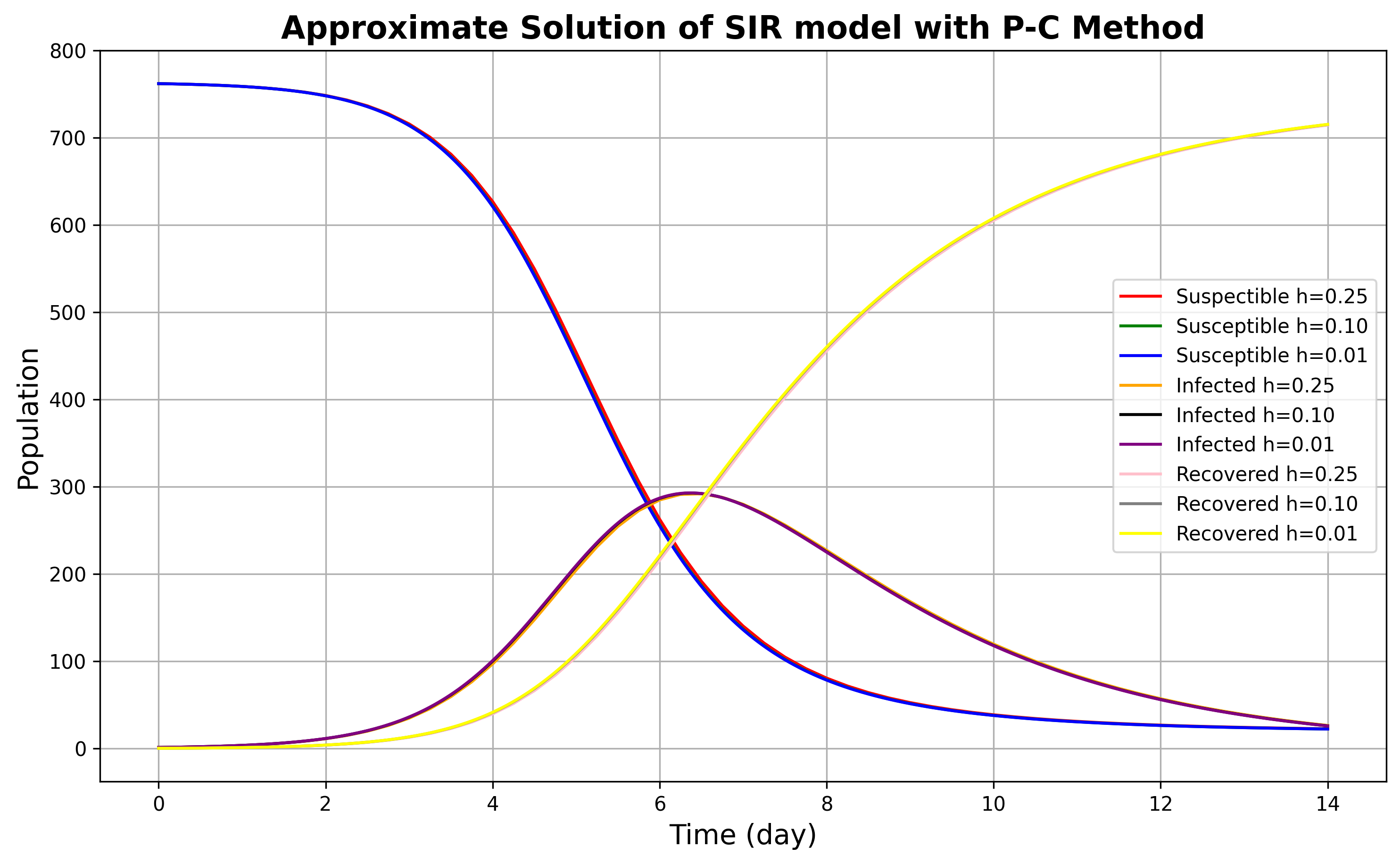}
        \caption{Solution of SIR model with P-C Method}
        \label{fig:Adams_Graph2}
    \end{subfigure}
    \caption{The graphs of the numerical solutions of the SIR model in Python, using Euler's method, RK4 method and the P-C method.}
    \label{fig:all_graphs}
\end{figure}

\begin{figure}[H]
    \centering
    \begin{subfigure}{0.3\textwidth}
        \centering
        \includegraphics[width=\linewidth]{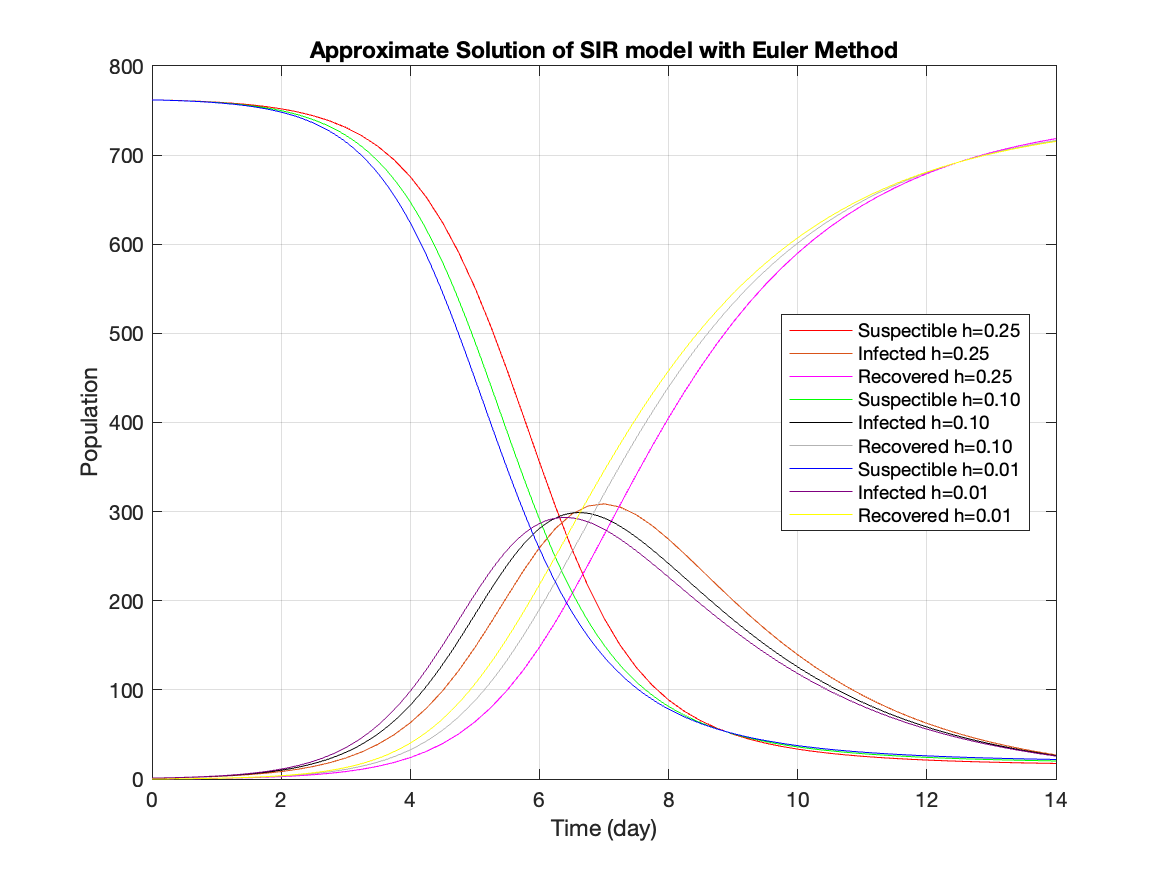}
        \caption{Solution of SIR model with Euler Method}
        \label{fig:euler_graph}
    \end{subfigure}
    \hfill 
    \begin{subfigure}{0.3\textwidth}
        \centering
        \includegraphics[width=\linewidth]{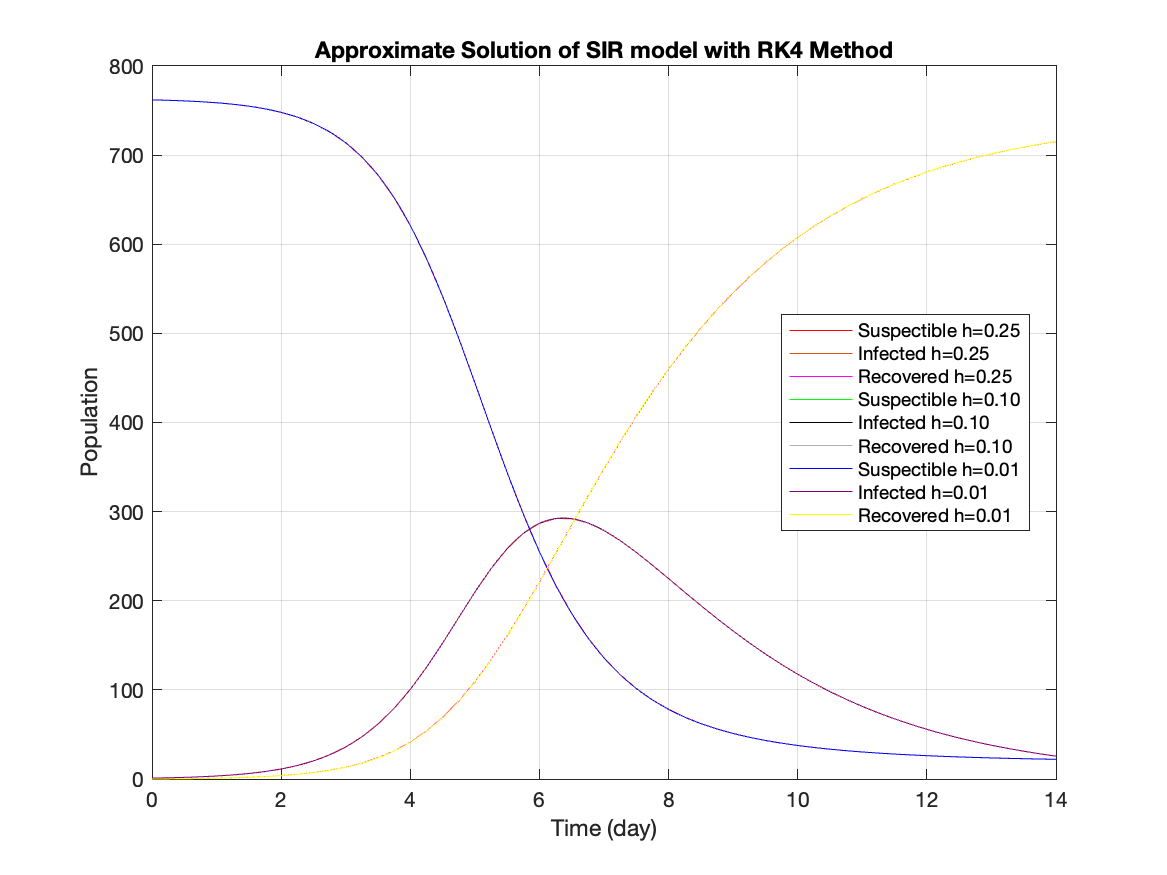}
        \caption{Solution of SIR model with RK4 Method}
        \label{fig:RK4_graph}
    \end{subfigure}
    \hfill 
    \begin{subfigure}{0.3\textwidth}
        \centering
        \includegraphics[width=\linewidth]{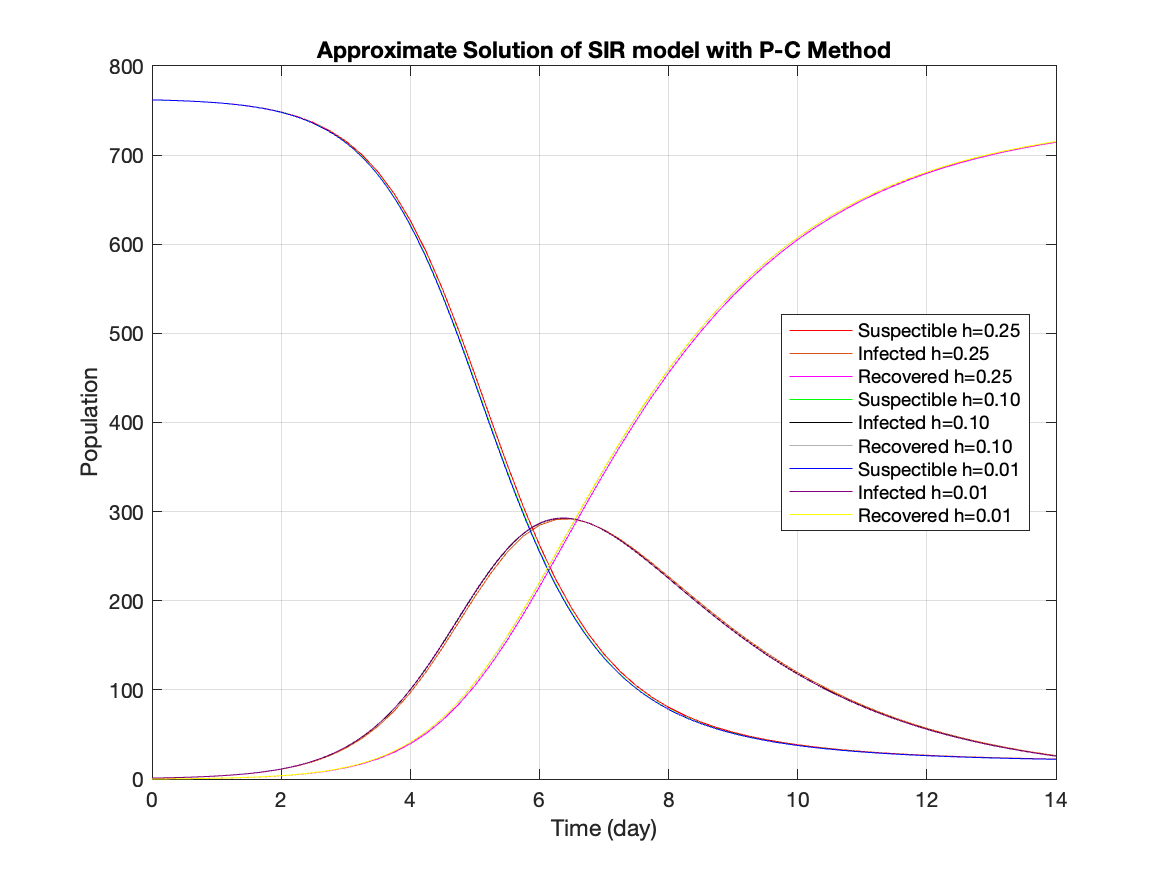}
        \caption{Solution of SIR model with P-C Method}
        \label{fig:Adams_Graph}
    \end{subfigure}
    \caption{The graphs of the numerical solutions of the SIR model in MATLAB, using Euler's method, RK4 method and the P-C method.}
    \label{fig:all_graphs}
\end{figure}

\begin{figure}[H]
    \centering
    \begin{subfigure}{0.3\textwidth}
        \centering
        \includegraphics[width=\linewidth]{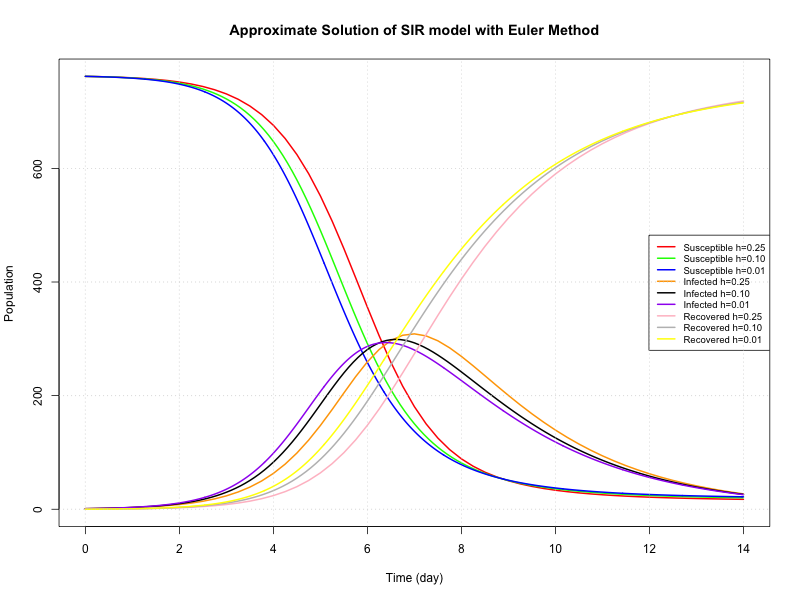}
        \caption{Solution of SIR model with Euler Method}
        \label{fig:euler_graph}
    \end{subfigure}
    \hfill 
    \begin{subfigure}{0.3\textwidth}
        \centering
        \includegraphics[width=\linewidth]{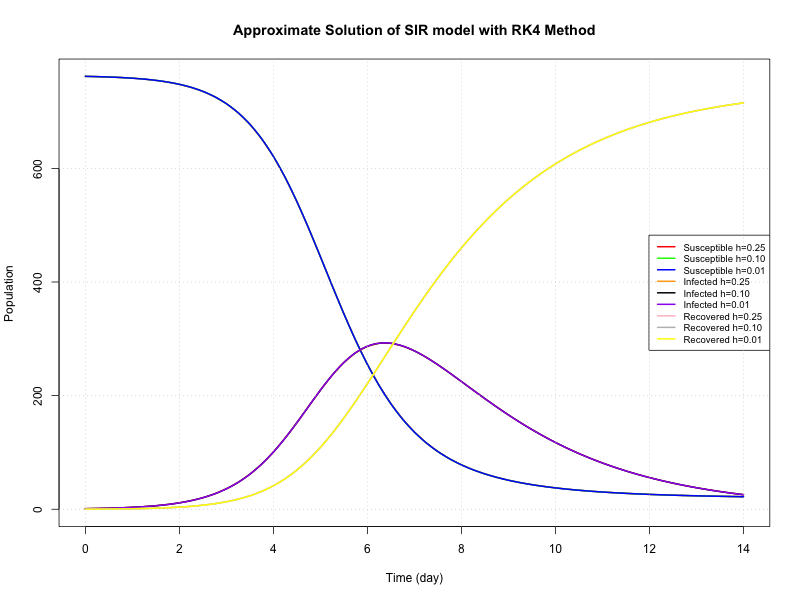}
        \caption{Solution of SIR model with RK4 Method}
        \label{fig:RK4_graph}
    \end{subfigure}
    \hfill
    \begin{subfigure}{0.3\textwidth}
        \centering
        \includegraphics[width=\linewidth]{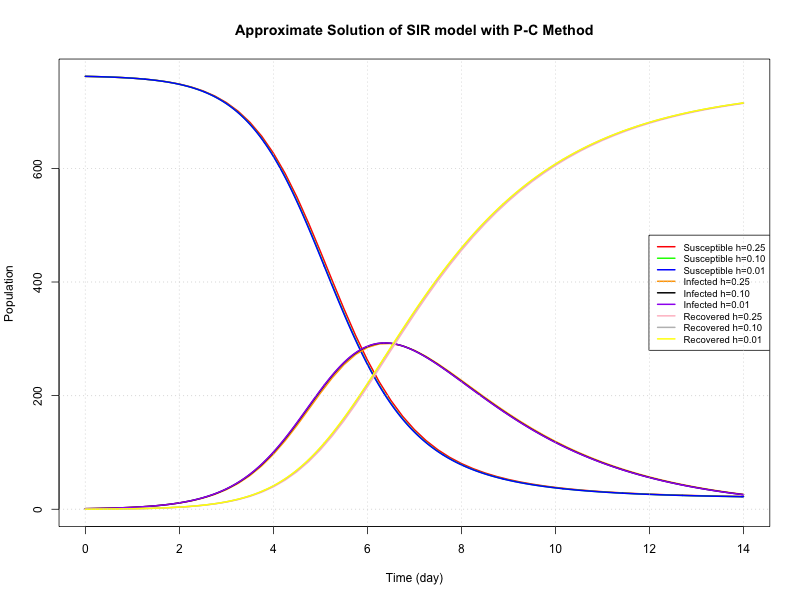}
        \caption{Solution of SIR model with P-C Method}
        \label{fig:Adams_Graph}
    \end{subfigure}
    \caption{The graphs of the numerical solutions of the SIR model in R, using Euler's method, RK4 method and the P-C method.}
    \label{fig:all_graphs}
\end{figure}

The errors of the RK4 method for the SIR model with respect to the ODE45 reference solution are given in Table 3.

\begin{table}[H]
    \centering
    \begin{tabular}{|>{\centering\arraybackslash}m{1.5cm}|c|c|}
        \hline
        \multirow{10}{*}{\rotatebox{90}{\textbf{RK4 Method}}} 
        &  & \textbf{MATLAB} \\ \cline{2-3}
        & S & h=0.25: $R^2$=0.9999998 \\
        &   & h=0.10: $R^2$=0.9999998 \\
        &   & h=0.01: $R^2$=0.9999998 \\ \cline{2-3}
        & I & h=0.25: $R^2$=0.9999996 \\
        &   & h=0.10: $R^2$=0.9999996 \\
        &   & h=0.01: $R^2$=0.9999996 \\ \cline{2-3}
        & R & h=0.25: $R^2$=0.9999998 \\
        &   & h=0.10: $R^2$=0.9999998 \\
        &   & h=0.01: $R^2$=0.9999998 \\ \hline
    \end{tabular}
    \caption{R$^2$ values of each numerical method in MATLAB (ODE45 as the reference solution and RK4) for the corresponding step sizes ($h$ values) for the SIR model. (Values are rounded to seven decimal places.)}
    \label{tab:numerical-errors}
    \end{table}

Finally, we present the run-times of the programs that solves the SIR model numerically.

\begin{table}[H]
\begin{center}
\begin{tabular}{|l|l|l|l|}
\hline
\textbf{} & \textbf{Python} & \textbf{MATLAB} & \textbf{R} \\
\hline
\textbf{Euler Method} 
& \begin{tabular}[c]{@{}l@{}}h=0.25: Time = 0.000089\\h=0.10: Time = 0.000125\\h=0.01: Time = 0.001185\end{tabular}
& \begin{tabular}[c]{@{}l@{}}h=0.25: Time = 0.008406\\h=0.10: Time = 0.006414\\h=0.01: Time = 0.015666\end{tabular}
& \begin{tabular}[c]{@{}l@{}}h=0.25: Time = 0.012420
\\h=0.10: Time = 0.007369
\\h=0.01: Time = 0.016969
\end{tabular} \\
\hline
\textbf{RK4 Method} 
& \begin{tabular}[c]{@{}l@{}}h=0.25: Time = 0.000244\\h=0.10: Time = 0.000555\\h=0.01: Time = 0.005373\end{tabular}
& \begin{tabular}[c]{@{}l@{}}h=0.25: Time = 0.015694\\h=0.10: Time = 0.024379\\h=0.01: Time = 0.030900\end{tabular}
& \begin{tabular}[c]{@{}l@{}}h=0.25: Time = 0.054937
\\h=0.10: Time = 0.031706
\\h=0.01: Time = 0.034045
\end{tabular} \\
\hline
\textbf{P - C Method} 
& \begin{tabular}[c]{@{}l@{}}h=0.25: Time = 0.000140\\h=0.10: Time = 0.000294\\h=0.01: Time = 0.002685\end{tabular}
& \begin{tabular}[c]{@{}l@{}}h=0.25: Time = 0.010834\\h=0.10: Time = 0.017722\\h=0.01: Time = 0.025391\end{tabular}
& \begin{tabular}[c]{@{}l@{}}h=0.25: Time = 0.020174
\\h=0.10: Time = 0.029155
\\h=0.01: Time = 0.018921
\end{tabular}
\\
\hline
\end{tabular}
\caption{Run-time of the software for each numerical method (seconds). (Values are rounded to six decimal places.)}
\end{center}
\end{table}

\section{Conclusions}
In this study, we compared the performance of three widely used programming environments, Python, MATLAB, and R, for the numerical solution of the SI and SIR epidemiological models given by (1) and (2), using Euler’s method, the RK4 method, and the P-C method. Our analysis was conducted in terms of both accuracy and computational efficiency.
For the SI model, which has an exact solution, we evaluated the accuracy of each numerical method using the $R^2$ value. Among the methods, RK4 consistently provided the highest accuracy across all software, achieving near-perfect results even for larger step sizes. The P-C method also provided highly accurate results, closely following RK4, while Euler’s method was the least accurate, especially for larger step sizes.
Regarding run-time performance, Python outperformed MATLAB and R in almost all cases, with significantly lower execution times, especially for finer step sizes. MATLAB showed moderate run-time performance, while R consistently exhibited the slowest execution among the three.
For the SIR model, where no exact analytical solution is available, the analysis focuses on both run-time performance and numerical behavior by employing a high-accuracy reference solution obtained via MATLAB’s ODE45 solver. The trends observed were consistent with those from the SI model: Python was the most efficient, while R was the least efficient. Additionally, as expected, the RK4 and P–C methods required more computational time than Euler’s method due to their more complex formulations. Furthermore, the SIR solutions computed using the RK4 method in MATLAB were compared against the ODE45 reference solution, and the results demonstrated that RK4 produces solutions that are very close to those obtained by ODE45, indicating that RK4 provides a reliable numerical approximation for the SIR model.
Overall, the study highlights that while all three software can be used effectively to solve SI and SIR models numerically, Python offers the best balance between speed and accuracy, especially when higher precision is required. These results obtained in this study can guide researchers in selecting appropriate computational tools for epidemic modeling based on their specific needs regarding accuracy, complexity, and performance.

\end{document}